\documentclass[a4paper,fullpage]{article}
\usepackage{amsmath}
\usepackage{amssymb}
\usepackage{amsthm}
\usepackage{times}
\usepackage{boxedminipage}
\usepackage[dvips]{graphicx}

\long\def\EqBox#1\BoxEq{\begin{equation}\fbox{$\displaystyle#1$}\end{equation}}
\long\def\EqRef#1#2\RefEq{$$#2\eqno{(#1)}$$}

\newtheorem{algo}{Algorithm}
 
\begin{document}
 
\begin{center}
\huge \textbf{ROBUST ROUTING AND CROSS-ENTROPY ESTIMATION}

\bigskip

H\'el\`ene Le Cadre 

\bigskip

\large{helene.lecadre@enst-bretagne.fr}

\bigskip
 
ENST, Bretagne\\

\bigskip

\large{September $9$, $2005$}
\end{center}

\section{Introduction} 

In this article we present a novel way to estimate the amounts of traffic on the Origin-Destination couples (OD couples). This new approach combines together a routing algorithm based 
on the principle of the shortest path and a recent technique of stochastic optimization called Cross-Entropy. The CE method was built at the origin, to tackle problems of rare-event simulation. However, its inventor, R. Rubinstein, realized soon that the underlying idea should be applied efficiently to combinatorial and multi-extremal optimization problems.\\
In a final part, we adapt a particular filtering algorithm in order to be able to dynamically estimate the evolution of the traffic on the OD couples.\\
The aim of this report is to highlight rather original ideas, however the choices of the prior distributions and some specific parameters may be quite arbitrary.

\section{A brief presentation of the CE method}

\subsection{Rare-Event Simulation}

Let $X\;=\;(X_1,...,X_N)$, be a random vector taking values in a space called $\mathcal{X}$. Let $\{f(.;v)\}$ be a family of parametric densities defined on the space $\mathcal{X}$, with respect to the Lebesgue measure. For any measurable function $H$, we can define:
$$\mathbf{E}[H(X)]\;=\;\int_{\Xi}\;H(x)\;f(x;v)\;dx\;.$$
The performance function will be called $S:\;\mathcal{X}\;\rightarrow\;\mathbb{R}$. For a fixed level $\gamma$, we are interessed in the probabilty of the event defined below:

$$l\;=\;\mathbf{P}_{u}(S(X)\;\geq\;\gamma)\;=\;\mathbf{E}_{u}[\mathbf{1}_{\{S(x)\;\geq\;\gamma\}}]\;.$$

If this probability is very small, for example not more than $10^{-5}$, the set $\{S(x)\;\geq\;\gamma\}$ will be called a \textit{rare event}. \\
A straightforward way to estimate $l$ may be to use crude Monte-Carlo simulation: let $(X_1,...,X_N)$ be a sample drawn from the density $f(.;u)$. Then, the estimator,

$$\frac{1}{N}\;\sum_{i=1}^{N}\;\mathbf{1}_{\{S(x)\;\geq\;\gamma\}}$$

is an unbiased estimator of $l$. However, if $\{S(x)\geq \gamma\}$  is a rare event, many indicator functions will remain equal to zero. As a result, we will be forced to simulate huge samples, which is rather costly and difficult to put in application.\\
Another way to get an estimate of $l$ might be to use importance sampling. We should draw $(X_1,...,X_N)$ from a density g , defined on the space $\mathcal{X}$. This density is nothing else than a mere change of measure. The estimator then becomes:

\begin{equation}
\label{is1}
\hat{l}\;=\;\frac{1}{N}\;\sum_{i=1}^{N}\;\mathbf{1}_{\{S(X_i)\;\geq\;\gamma\}}\;\frac{f(X_i;u)}{g(X_i)}\;.
\end{equation}

The optimal density g, is defined by:

\begin{equation}
\label{is2}
g^{\star}(x)\;=\;\frac{\mathbf{1}_{\{S(x)\;\geq\;\gamma\}}\;f(x;u)}{l}\;.
\end{equation}

Substituting $(\ref{is2})$ in $(\ref{is1})$, we get:

$$\mathbf{1}_{\{S(X_i)\;\geq\;\gamma\}}\;\frac{f(X_i;u)}{g^{\star}(X_i)}\;=\;l\;,\;\forall\;i\;.$$

But, $l$ is a constant. As a result, the estimator defined in $(\ref{is1})$ has zero variance. Nevertheless, $g^{\star}$ depends on the unknown parameter $l$. The idea is in fact, to choose g in a parametric family of densities $\{f(.;v)\}$. The problem is now to determine the optimal parameter $v$, such that the distance between $g^{\star}$ and $f(.;v)$ should be minimized.\\
A well-known "distance" between two densities g and h, is the Kullback-Leibler "distance":

\begin{equation}
\mathcal{D}(g,h)\;=\;\mathbf{E}_{g}[ln\;\frac{g(X)}{h(X)}]\;=\;\int\;g(x)\;ln\;g(x)\;dx\;-\;\int\;g(x)\;ln\;h(x)\;dx\;.
\end{equation}
 
Minimizing the Kullback-Leibler distance between $g^{\star}$ and $f(.;v)$, is equivalent to solve the following problem:

\begin{equation}
\label{Kull}
\displaystyle{\arg \max \limits_{v}}\;\int\;g^{\star}(x)\;ln\;f(x;v)\;dx\;.
\end{equation}

Substituting $(\ref{is2})$ in $(\ref{Kull})$, we get:
\begin{equation}
\displaystyle{\arg \max \limits_{v}}\;D(v)\;=\;\displaystyle{\arg \max \limits_{v}}\;\mathbf{E}_{u}[\mathbf{1}_{\{S(X)\;\le\;\gamma\}}\;ln\;f(X;v)]\;.
\end{equation}

But, in fact we can estimate $v^{\star}$, by solving the following stochastic program:

\begin{equation}
\label{pbfond}
\displaystyle{\arg \max \limits_{v}}\;\hat{D}(v)\;=\;\displaystyle{\arg \max \limits_{v}}\;\frac{1}{N}\;\sum_{i=1}^{N}\;[\mathbf{1}_{\{S(X_i)\;\geq\;\gamma\}}\;ln\;f(X_i;v)]\;.
\end{equation}

If $\hat{D}$ is convex and differentiable in $v$, we just have to solve the following problem:

\begin{equation}
\frac{1}{N}\;\sum_{i=1}^{N}\;[\mathbf{1}_{\{S(X_i)\;\ge\;\gamma\}}\;\nabla_{v}\;ln\;f(X_i;v)]\;=\;0\;.
\end{equation}

The solution can often be calculated analytically, which is one of the great advantages of this approach.

\subsection{Application of the CE to optimization}

Usually, in the field of optimization we try to solve problems of the form:
\begin{equation}
\label{opt}
S(x^{\star})\;=\;\gamma^{\star}\;=\;\displaystyle{\arg \max \limits_{x\in\mathcal{X}}}\;S(x)\;.
\end{equation}

The genius of the Cross-Entropy method lies in the fact that it is possible to associate with each optimization problem of the form $(\ref{opt})$, a problem of estimation, called \textit{associated stochastic problem (ASP)}. We will start by defining a collection of indicator functions $\{\mathbf{1}_{\{S(x)\;\le\;\gamma\}}\}_{\gamma\;\in\;\mathbb{R}},$ on the space $\mathcal{X}$. Then, we will define a parametric family of densities $\{f(.;v),\;v\;\in\;\mathcal{V}\}$ on the space $\mathcal{X}$. Let $u\;\in\;\mathcal{V}$. We will associate with $(\ref{opt})$, the following stochastic estimation problem:

\begin{equation}
\label{ASP}
l(\gamma)\;=\;\mathbf{P}_{u}(S(X)\;\geq\;\gamma)\;=\;\sum_{x}\;\mathbf{1}_{\{S(x)\;\geq\;\gamma\}}\;f(x;u)\;=\;\mathbf{E}_{u}[\mathbf{1}_{\{S(x)\;\geq\;\gamma\}}]\;,
\end{equation}

If $\gamma\;=\;\gamma^{\star}$, a natural estimator of the reference parameter $v^{\star}$  is: 
\begin{equation}
\hat{v^{\star}}\;=\;\displaystyle{\arg \max \limits_{v}}\;\frac{1}{N}\;\sum_{i=1}^{N}\;\mathbf{1}_{\{S(x)\;\geq\;\gamma\}}\;ln\;f(X_i;v)\;,
\end{equation}

where the $X_i$ are drawn from the density $f(.,u)$. If $\gamma$ is very close to $\gamma^{\star}$, then $f(.;v^{\star})$ assigns most of its probability mass close to $x^{\star}$. In fact, in this case, we will have to choose $u$ so that $\mathbf{P}_{u}(S(X)\;\geq\;\gamma)$ is not too small. We can infer that $u$ and $\gamma$ are closely linked. \\
We will use a two level procedure. Indeed, we will construct two sequences $\hat{\gamma_1},...,\hat{\gamma_T}$ and $\hat{v_0},\hat{v_1},...,\hat{v_T}$ such that $\hat{\gamma_T}$ is close to $\gamma^{\star}$ and $\hat{v_T}$ is such that the density assigns most of its mass in the state which maximizes the performance. \\
The algorithm follows a two-step strategy:

\begin{algo}
\begin{itemize}
\item Define $\hat{v_0}\;=\;u$. Set $t\;=\;1$.
\item Generate a sample $(X_1,...,X_N)\;\sim\;f(.;v_{t-1})$. Compute the $(1-\rho)$-quantile of the performance S, which can be estimated by: $$\hat{\gamma_t}\;=\;S_{[(1-\rho)N]}\;.$$
\item Use the same sample $(X_1,...,X_N)$, to solve $(\ref{pbfond})$. Call the solution $v_t$.
\item If for some $t\;\geq\;d$, d fixed, $$\hat{\gamma_t}=...=\hat{\gamma_{t-d}}\;,$$ \textbf{Stop}; otherwise set $t=t+1$ and reiterate from Step $2$.
\end{itemize}
\end{algo}

\section{The Model}

The network we study is composed of $p$ nodes and $n$ arcs. At first, we will suppose that there exists an arc between each couple of nodes. Furthermore, we make the difference between the two couples $(i,j)$ and $(j,i),\;\forall\;i,j\;\in\;\{1,...,n\},\;i\;\neq\;j.$. Consequently, $$n\;=\;p^{2}\;-\;p\;.$$
Recall that an arc is a directed link. In the rest of this article, we will make the hypothesis that the network is directed. 

\bigskip

Our work could be separated into two different parts. Firstly, we have to deal with a simulation part. In this section, we will choose an initial vector of the amounts of traffic on the OD couples. Our aim will be to minimize the global sum of the costs, which are associated to each arc of the network. Secondly, we will have to cope with an estimation part. Indeed, we will suppose that the initial costs are equal to those obtained by the simulation part. The idea then, will be to find the estimator $\hat{X}(t)$ from which we could infer an arc estimator $\hat{Y}(t)$ minimizing the distance to the vector $Y(t)$, obtained in the simulation part. This raises the crucial question of the identifiability of the vector $X(t)$. Is there unicity of the associated vector $\hat{X}(t)$ or, is it only an element in a vast variety?

\subsection{Simulation}

In this part, we associate a cost function to the network, which means that we give a cost to each arc of the network. This cost function is drawn from a parametric families of densities. Which means that we have to determine the optimal parameter of the density function. What's more, this cost may or, may not be, proportionnal to the amount of traffic on each arc. However, it is more unconventional to suppose that the costs depend on the arc traffic. Let $Y(t)$ be the vector which contains the amount of traffic on each arc of the network and $C(t)$ the vector which represents the costs associated with each arc. 
In the most general case, we have:
\begin{equation}
C(t)\;=\;F(Y(t))\;,
\end{equation}
where F is supposed to be continuous and differentiable. 

\bigskip

\begin{center}
\begin{figure}[h]
\includegraphics[scale=0.5]{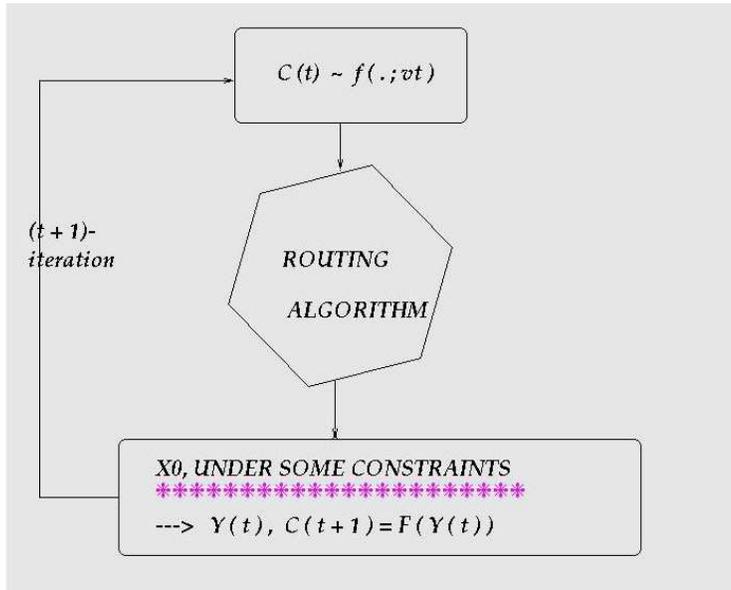}
\caption{Simulation.}\label{simu}
\end{figure}
\end{center}
\bigskip

Then, using these costs, we will use a routing algorithm based on the principle of the shortest path in order to find the shortest paths between each OD couple. A path is represented by the nodes or the arcs which it is made of.\\
The principle  of the routing algorithm we use is rather simple. The first observation to make is that every subpath of a shortest path is necessarily itself a shortest path. \\
Let $e_{i,j}$ be the arc linking the nodes i and j. If the path composed of the arcs $\{e_{i,j},\;e_{j,k},\;e_{k,l},...,e_{p,q}\}$ is the shortest path linking node $i$ to node $q$. Then, $e_{i,j}$  must be the shortest path between $i$ and $j$. $e_{j,k}$ must be the shortest path between $j$ and $k$, and so on... \\
$e_{i,j}$ is called \textit{basic arc}, iff it is the shortest path between $i$ and $j$. Consequently, each shortest path must be composed exclusively of basic arcs. The aim of our routing algorithm is to substitute to each arc which is not basic, a basic arc. Let $d_{i,j}\;=\;C(t;i,j)$ be the distance or weight associated to the arc 
$\{i,j\}$, at the instant $t$. Let $j$ be the indice of a node of the network, then:

$$\forall\;i\;,k\;\in\;\textrm{Network}-\{j\}\;,\;d_{i,k}\;\leftarrow\;min\{d_{i,k},d_{i,j}+d_{j,k}\}\;.$$

The algorithm tests all the couple of nodes $(i,k)$, which are neighbors of $j$, while $j$ takes each node of the network as its own value. 

\bigskip

\begin{algo}
\textbf{Input}: $C(t)$, vector of the costs.
\begin{itemize}
\item If $d_{i,k}\;\geq\;d_{i,j}+d_{j,k}$, do not change anything.
\item If $d_{i,k}\;\leq\;d_{i,j}+d_{j,k}$, create an arc linking $i$ to $k$ and associate the weight $d_{i,k}\;=\;d_{i,j}+d_{j,k}$.
\end{itemize}
\textbf{Output}: the shortest paths between each OD couple.
\end{algo}

The algorithm also gives us the shortest distances associated to each OD couple. But, these distances are only rought estimators of the amounts of traffic on each OD couple. Indeed, more than one link, can be shared by different shortest paths linking different OD couples. As a result, the total amount of traffic generated by \textbf{one} OD couple usually represents only a fraction of the total traffic flowing through the arcs which composed the path.

\bigskip

\begin{figure}[h]
\includegraphics[scale=0.5]{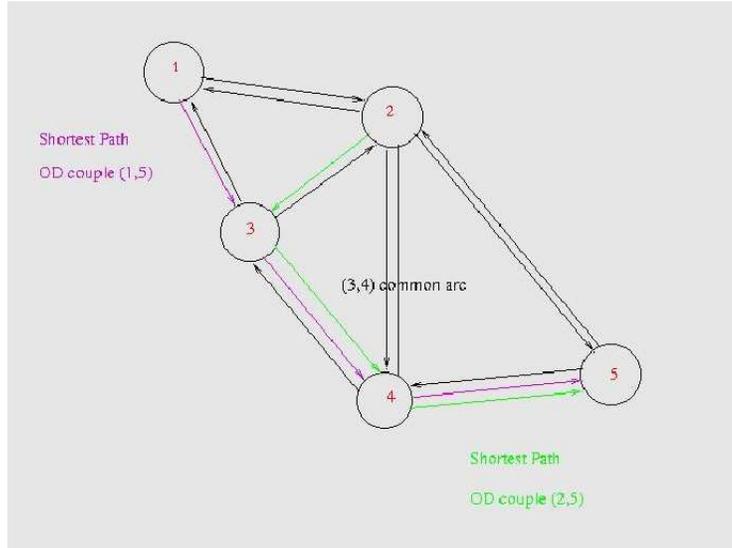}\hfill
\caption{Connection between OD couples and arcs.}\label{est}
\end{figure}

In order to solve this crucial problem, we will associate to the vector which contains the volumes of traffic on the OD couples, called $X_0(t)$, an estimator of the volumes of traffic flowing through the arcs, which we will note $Y(t)$. Indeed, if we use the routing algorithm, it is quite easy to deduce $Y(t)$ from $X_0(t)$. Our goal will be to solve the following optimization problem:

\begin{equation}
min\;\sum_{i=1}^{n}\;|C_i(t)|\;.
\end{equation}

\subsection{Estimation}

In this part, the performance function is defined by:
\begin{equation}
S(X(t))\;=\;\displaystyle{\arg \max \limits_{\hat{X(t)}}}\;\frac{1}{\|Y(t)-\hat{Y}(t)\|_{2}}\;.
\end{equation}

\bigskip

\begin{center}
\begin{figure}[h]
\includegraphics[scale=0.5]{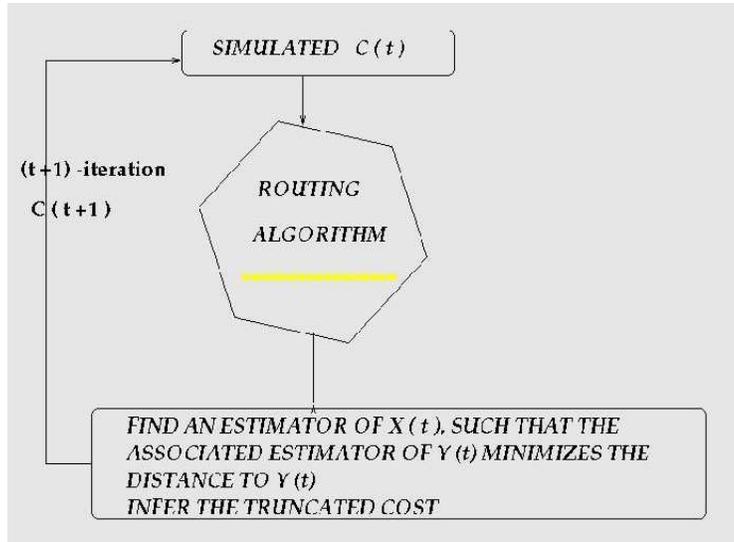}\hfill
\caption{Estimation}\label{est}
\end{figure}
\end{center}

\bigskip

The first observation is that $S$ is an implicit function of $X(t)$, e.g. we can't get any exact analytical expression of $S$. Consequently, we will have to resort to use simulation.

\bigskip

What's more, we suppose that the vector $X(t)$ is generated from an exponential density whose parameter is totally unknown. 

\begin{equation}
X(t)\;\sim\;\mathcal{E}(\lambda),\;\lambda\;\in\;\mathbb{R}^{n}_{+}\;.
\end{equation}

We apply the Cross-Entropy method to our problem. At each iteration, we generate a random sample $(X^{(1)},...,X^{(N)})\;\sim\;\mathcal{E}(\lambda),\;\lambda\;\in\;\mathbb{R}^{n}_{+}\;.$ 

\bigskip

\textbf{Hypothesis}: Each component of $X(t)$, which represents the amount of traffic on an OD couple, will be supposed to be independent of the others. 

\bigskip

Practically, each random vector from the sample will be stocked in a big matrix.

\begin{equation}
X(t)\;=\;\left[
\begin{array}{c c c}
X_1^{(1)} & \cdots & X_1^{(N)}\\
X_2^{(1)} & \cdots & X_2^{(N)}\\
\vdots & \vdots & \vdots\\
X_n^{(1)} & \cdots & X_n^{(N)}\\
\end{array}
\right].
\end{equation}

The joined densities of the vectors $X^{(j)},\;j=1,...,N$, are typically of the form:
\begin{eqnarray}
f(X^{(j)},\lambda)\;&=&\;\prod_{i=1}^{n}\;\lambda_i\;\exp[-\lambda_i\;X_{i}^{(j)}]\;\mathbf{1}_{\mathbb{R}^{+}}(X_{i}^{(j)})\;,\nonumber\\
{} &=& \;\prod_{i=1}^{n}\;\lambda_i\;\exp[-\lambda_i\;X_{i}^{(j)}]\;\mathbf{1}_{\{min(X_{i}^{(j)})\;\geq\;0\}}\;.
\end{eqnarray}

As a result, we will have to solve the following problem:

\begin{equation}
\frac{1}{N}\;\sum_{i=1}^{N}\;\mathbf{1}_{\{S(X^{(i)})\;\geq\;\hat{\gamma_t}\}}\;\nabla_{\lambda}\;ln\;f(X^{(i)};\lambda)\;=\;0\;.
\end{equation}

After some computations, we get:

\begin{equation}
\lambda_j\;=\;\frac{\sum_{i=1}^{N}\mathbf{1}_{\{S(X^{(i)})\;\geq\;\hat{\gamma_t}\}}}{[\sum_{i=1}^{N}\mathbf{1}_{\{S(X^{(i)})\;\geq\;\hat{\gamma_t}\}}\;X^{(i)}_j]},\;\forall\;j=1,...,n\;.
\end{equation}

\bigskip

\textbf{Remark}.\\
If we generate the random vectors $X_{i},\;i=1,...,N$  from a truncated exponential, we can give some maximal boundaries on the OD volumes of traffic.\\
Recall that a truncated exponential is of the form:
$$f(x;\lambda,b)\;=\;\frac{\lambda\;\exp[-\lambda\;x]}{1-\exp(-\lambda\;b)}\;\mathbf{1}_{[0,b]}(x)\;.$$

Under this assumption, we have to cope with the following system:

\begin{equation}
\frac{\displaystyle{\sum_{i=1}^{N}}\;\mathbf{1}_{\{S(X^{(i)})\;\geq\;\hat{\gamma_t}\}\;X_j^{(i)}}}{\displaystyle{\sum_{i=1}^{N}}\;\mathbf{1}_{\{S(X^{(i)})\;\geq\;\hat{\gamma_t}\}}}\;-\;\frac{1}{\lambda_j}\;+\;\frac{b_j}{\exp(\lambda_j\;b_j)-1}\;=\;0\;,\;\forall\;j\;\in\;\{1,...,n\}\;.
\end{equation}

This system is non-linear, that's why we use the well-known iterative Newton's method to solve it.

\bigskip

\begin{center}
\begin{figure}[h]
\includegraphics[scale=0.7]{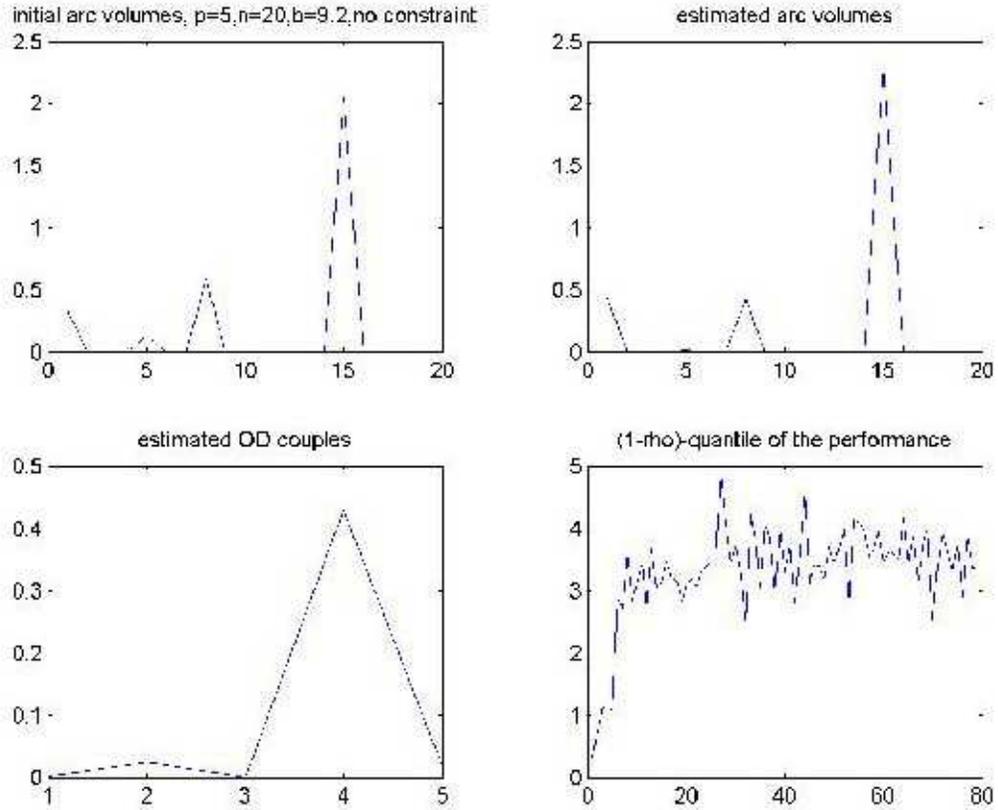}
\caption{Simulation, $5$ nodes network, $20$ arcs.}\label{simu5}
\end{figure}
\end{center}

\bigskip

\begin{center}
\begin{figure}[h]
\includegraphics[scale=0.7]{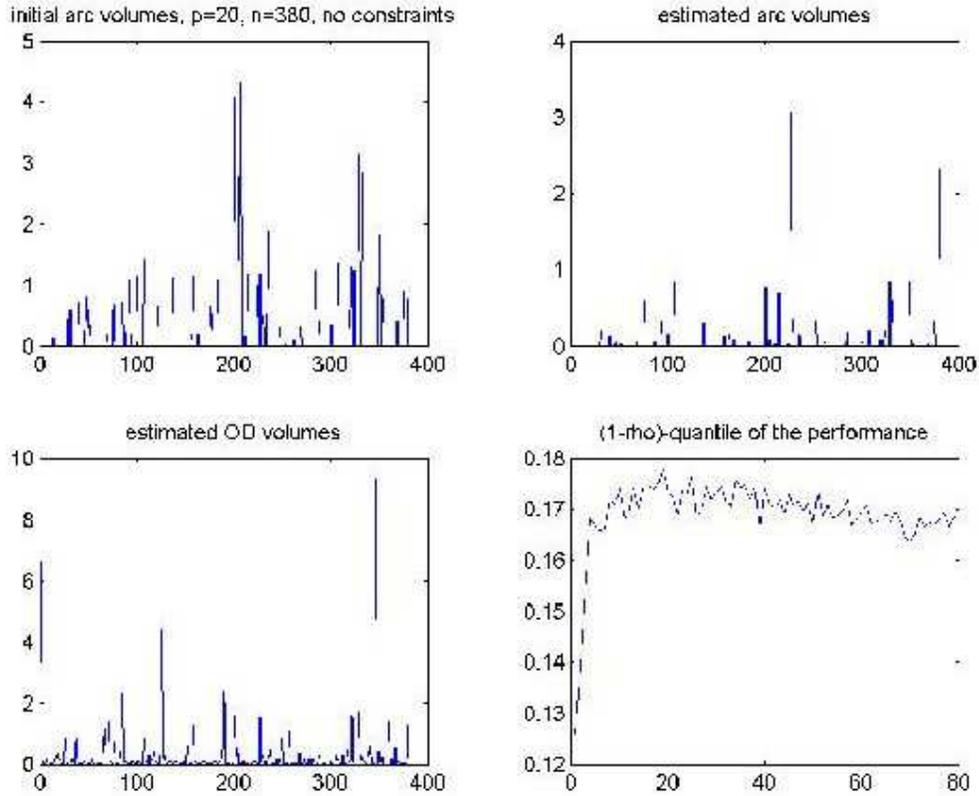}
\caption{Simulation, $20$ nodes network, $380$ arcs.}\label{simu20}
\end{figure}
\end{center}

\section{The problem of Identifiability}

We observe that the amount of traffic which flows through each arc of the network is equal to the sum of the amounts of traffic flowing on each OD couple which owns this arc in its shortest path. Remember that, thanks to the routing algorithm, we associate to each OD couple, a unique shortest path. 
Mathematically, we can express this remark under the following expression:

\begin{equation}
\sum_{j=1}^{n}\mathbf{1}_{\{e_i\;\in\;\textrm{OD couple number j's path}(t)\}}\;X_j(t)\;=\;Y_i(t)\;,\;\forall\;i\;\in\;\{1,...,n\}.
\end{equation}

More generally, we get:

\begin{equation}
\label{routage}
Y(t)\;=\;A(X(t))\;X(t)\;.
\end{equation}

Indeed, in the most general case, the routing matrix A relies on the volumes of traffic flowing through each arc at the instant $t$. But, these arc volumes rely themselves on the OD volumes, $X(t)$. As a result, the routing matrix $A(t)$, is a function of $X(t)$.\\ 
A is uniquely made of binary elements: $0$ and $1$. More explicitely, $A(X(t);i,j)\;=\;1$ iff, the arc numbered $i$ belongs to the shortest path associated to the OD couple numbered $j$, at time $t$.\\
What's more, the routing algorithm do not use every arc. Consequently, many rows of the routing matrix equal zero. The associated components in the arc vector $Y(t)$ are at the same time, null.\\
But, if we suppress the zero rows of $A(X(t))$ and the zero components of $Y(t)$, this leads us to solve a rectangular system of equations. This system is under-determined, that's why we can't guarantee the existence of a unique solution.

\bigskip

We can conclude that there is no identifiability between the arc volumes $Y(t)$, and the OD volumes $X(t)$, at a given time. Indeed, if we take a fixed $X(t)$, we get a unique associated $Y(t)$, since the routing algorithm determines a unique shortest path between each OD couples. Reciprocally, if we take some fixed arc volumes, $Y(t)$, we can't guarantee the unicity of the solutions of $(\ref{routage})$. That's why, we can't assert that the associated $X(t)$ is perfectly unique.

\bigskip

A good idea to tackle this problem, is to suppose that some of the OD couples do not accept any traffic. That is, that they remain equal to zero. The goal is to reduce the number of positive OD couples so as to get a system whose routing matrix $A(t)$, is square or not too far.

\bigskip

The first approach is to suppose that some pre-determined OD couples are exculded.

To begin, we may partition, a little arbitrarily, the set of the OD couples into two parts. In the first one, lie the OD couples which remain always equal to zero. And, in the second part, 
we will suppose that there is some traffic flowing through these OD couples.\\

We need to generate a sample $Z$. The components of $Z$ are independent of each other and generated from a Bernoulli density whose parameter is pre-determined.
Then, if $Z(i)=0$, the OD couple number $i$, do not accept any traffic. 

\bigskip
\begin{center}
\begin{figure}[h]
\includegraphics[scale=0.7]{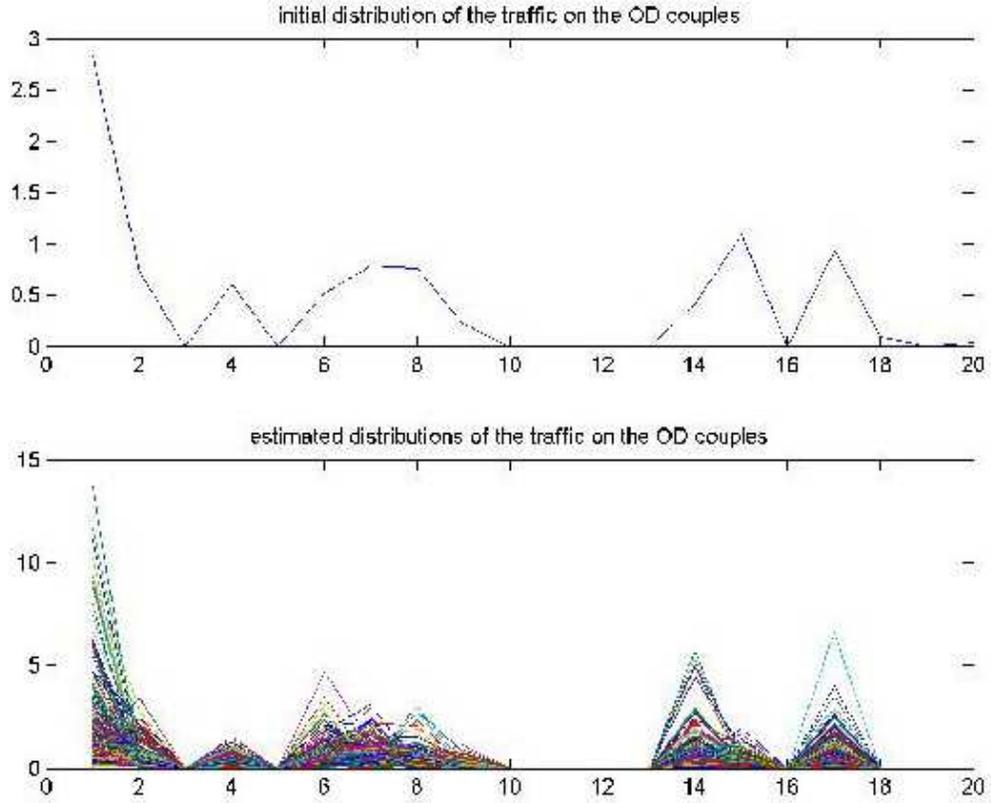}
\caption{Introduction of constraints on the OD couples: the zero OD couples are pre-determined. Identifiability of X(t) with $\frac{2}{3}$ of the OD couples used.}\label{contrmeth1}
\end{figure}
\end{center}

\bigskip

The time required to perform this simulation is of about $1$ minute. 

\bigskip 

A second point of vue should be to suppose that we know at the beginning that only $K$ OD couples, $K\;\geq\;n$, are positive.  
So, we need to modify our CE algorithm. We need to introduce a matrix, $Z(t)$:

\begin{equation}
Z(t)\;=\;\left[
\begin{array}{c c c}
Z_1^{(1)} & \cdots & Z_1^{(N)}\\
Z_2^{(1)} & \cdots & Z_2^{(N)}\\
\vdots & \vdots & \vdots\\
Z_n^{(1)} & \cdots & Z_n^{(N)}\\
\end{array}
\right].
\end{equation}

To be more explicit, the $i^{th}$ row of $Z$ is generated from $\mathcal{B}(p_{i,j}(t)),\;j\;\in\;\{0,1\},\;\forall\;i\;\in\;\{1,...,n\}\;$.  
In fact, each row is generated independently from a Bernoulli density whose parameter is specific, conditional upon the fact that $\sum_{j=1}^{n}X_j^{(i)}\;=\;K,\;i\;\in\;\{1,...,N\}\;.$
K, is a fixed number. It may be as we have already stated, a certain propportion of OD couples, but it may also take into account some other constraints. \\
The first idea to deal with such a constraint is to generate a random vector $X_1^{(i)},...,X_n^{(i)}$. Each component are drawn independently from a Bernoulli density. The sample is accepted iff, $\sum_{j=1}^{n}X_j^{(i)}\;=\;K,\;i\;\in\;\{1,...,N\}\;.$
However, when $n$ becomes higher than $10$, it takes a prohibitive time! In fact, the best solution is to generate independant Bernoulli random variables from $\mathcal{B}(p_1(t)),\mathcal{B}(p_2(t)),...$, respectively, until $K$ unities or $n-K$ zeros are generated. Then, the remaining elements are put equal to zero or one, respectively.
However, the updating formula for the parameters of the Bernoulli densities remain exactly of the form:

\begin{equation}
p_{i,1}(t)\;=\;\frac{\sum_{k=1}^{N}\;\mathbf{1}_{\{S(X^{(k)}\;\geq\;\hat{\gamma_t}\}}\;\mathbf{1}_{\{X_i^{(k)}\;=\;1\}}}{\sum_{k=1}^{N}\;\mathbf{1}_{\{S(X^{(k)}\;\geq\;\hat{\gamma_t}\}}}\;,\;\forall\;i\;\in\;\{1,...,n\}\;.
\end{equation}

In fact, now, $X(t)$ and $Z(t)$ are closely linked. Indeed, if $Z_i^{(j)}\;=\;0$ then, $X_i^{(j)}\;=\;0$, which means than there is no traffic on the OD couple number $i$ for the $j^{th}$-sample. 

\bigskip
\begin{center}
\begin{figure}[h]
\includegraphics[scale=0.7]{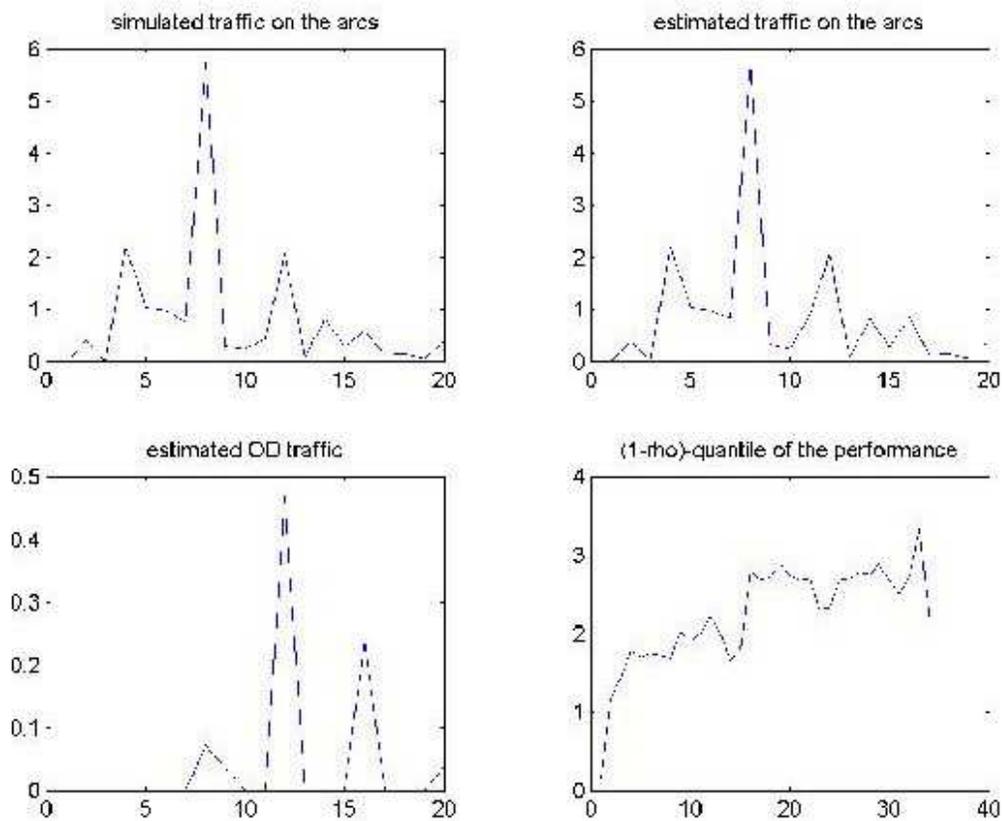}
\caption{Simulation, the number of zero OD couples is pre-determined, $K=\frac{2}{3}.$}\label{simumeth2}
\end{figure}
\end{center}

\bigskip

\begin{center}
\begin{figure}[h]
\includegraphics[scale=0.7]{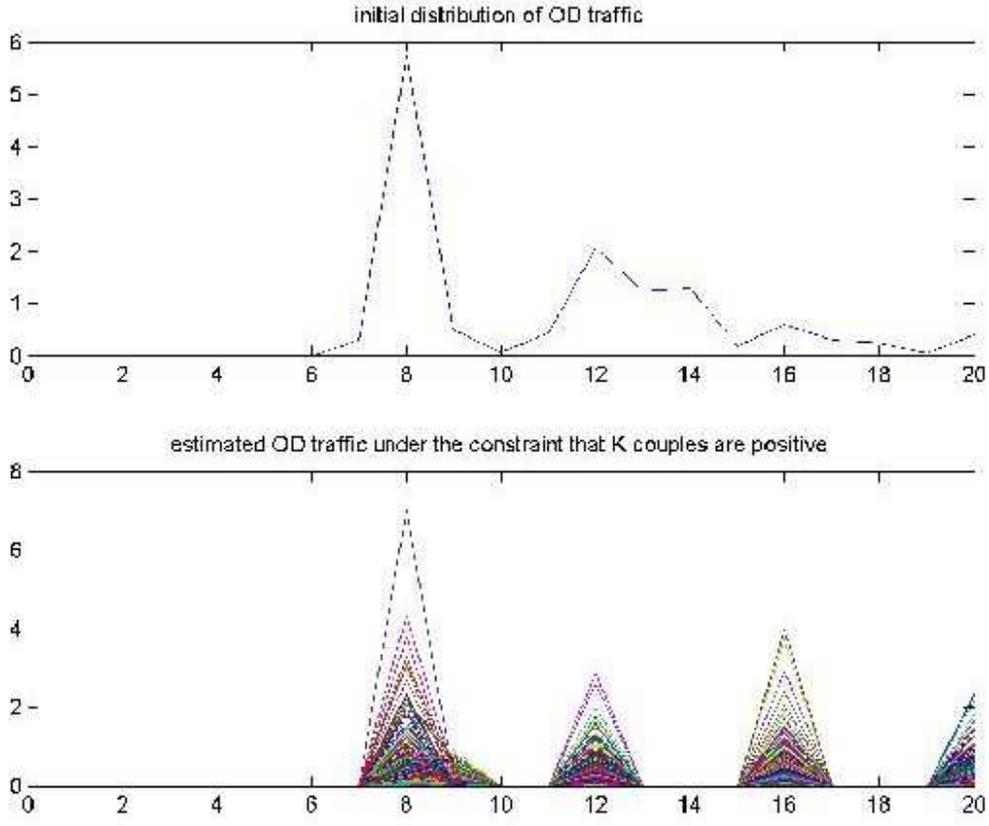}
\caption{Identifiability of X(t), $K=\frac{2}{3}$.}\label{simumeth2b}
\end{figure}
\end{center}
The time required to perform this simulation is of about $2$ minutes.

\section{Dynamic estimation}

We have previously determine an estimator of the amounts of traffic flowing through the OD couples at the specific instant t. We should ask ourselves whether it is possible to determine the trajectories associated to the vector X(t). Particle filtering appear to us to be an interesting approach. \\

\subsection{Presentation of Particle filtering}

Particle filtering is a well-known technique based on sequential Monte-Carlo approach. It is a technique for implementing a recursive bayesian filter by Monte-Carlo simulations. The key ideea is to represent the required posterior density function by a set of random samples with associated weights and to compute estimates based on these samples and weights.\\

We will generate a random measure $\{C_{0:k}^{i},w_{k}^{i}\}_{i=1,...,M}$ that characterises the posterior pdf $p(C_{0:k}\;|\;\hat{X}_{1:k})$. $\{C_{0:k}^{i},i=0,...,M\}$ is a set of vectors with associated weights $\{w_{k}^{i},i=1,...,M\}$ (the weights are themselves vectors of weights). $C_{0:k}=\{C_{j},j=0,...,t\}$ is the set of all states up to time t. The weights are normalised such that, $$\sum_{i=1}^{M}w_{k}^{i,j}\;=\;1,\;\forall\;j\;\in\;\{1,...,n\}.$$

Then, the posterior density at t can be approximated as:

\begin{equation}
\label{appro}
p(C_{0:k}\;|\;\hat{X}_{1:k})\;\thickapprox\;\sum_{i=1}^{M}\;w_k^{i}\;\delta(C_{0:k}-C_{0:k}^{i})\;.
\end{equation}

The weights are chosen using the principle of Importance Sampling. Let $C^{i}\;\sim\;q(.),\;i=1,...,M$ be samples generated from a proposal $q(.)$, called Importance sampling density.
By successive approximations, it is shown in $[7]$ that the weights are recursively obtained by the following formula:

\begin{equation}
w_{k}^{i}\;\propto\;w_{k-1}^{i}\;\frac{p(\hat{X_k}|C_{k}^{i})p(C_{k}^{i}|C_{k-1}^{i})}{q(C_{k}^{i}|C_{k-1}^{i},\hat{X_k})}\;.
\end{equation}

It can be shown that as $M\;\rightarrow\;\infty$, the approximation $(\ref{appro})$ approaches the true posterior density $p(C_{k}\;|\;\hat{X}_{1:k})$. 

However, there is a major drawback to use particle filtering techniques. Indeed, a common problem is the degeneracy problem. After a few iterations, all but one particle will have negligable weight. It has been shown that the variance of the importance weights can only increase over time, and thus it is impossible to avoid the degeneracy phenomenon. This degeneracy implies that a large computational effort is devoted to updating particles whose contribution to the approximation to $p(C_{k}|\hat{X_{1:k}})$ is almost zero.\\
Consequently, we have to use resampling mechanisms. The basic idea behind resampling is to eliminate particles which have small weights and to concentrate on particles with large weights.\\

\subsection{State model and observation equation}

The traffic flow will be modelled as a stochastic hybrid system with discrete states. The observation equation is rather simple to get. Indeed, we have:

\begin{equation}
\hat{X_{k}}\;=\;\Xi(C(k))\;.
\end{equation}

Where, $\Xi$ is a quite complex function which represent the whole algorithm. \\
The difficulty now, is to build a state model. Suppose the flow can be decomposed in small entities (for example packets). 
We note: $\{Q_{l',k}| l'\;\in\;\{\textrm{arcs of the network}\},\;i\rightarrow l'\}$, the number of packets going out of the arc i, during the time interval $[t_k,t_{k+1}[$. $\{Q_{l,k}| l\;\in\;\{\textrm{arcs of the network}\},\;l\rightarrow i\}$, is the number of packets arriving on the arc i on $[t_k,t_{k+1}[$. 

The conservation of the flow lets us write:

\begin{equation}
\label{conserve}
Y_{i}(k+1)\;=\;Y_{i}(k)\;+\;\sum_{\{arcs\;l|l\rightarrow i\}}Q_{l,k}\;-\;\sum_{\{arcs\;l'|i\rightarrow l'\}}Q_{l',k}\;.
\end{equation}

In fact,

\begin{equation}
Q_{i,k}\;=\;min(S_{i,k};R_{i,k+1})\;.
\end{equation}

$S_{i,k}$ is called \textit{sending function}. It expresses how many among the $Y_{i}(k)$ packets in the arc i at k are at a distance less than a fixed boundary called $\beta$. Suppose the interaction between the packets is negligible and their location is uniformly distributed over the arc. $S_{i,k}$ is then a random binomial variable with $Y_{i}(k)$ drawings, with probability of success $\frac{\beta}{\textrm{arc i length}}$, or an approximation, since we don't know exactly the length of the arc number i. 

\bigskip

The \textit{receiving function} is defined by:

\begin{equation}
R_{i,k+1}\;=\;\sum_{\{arcs\;l|i\rightarrow l\}}[Y_{l}^{max}(k)\;+\;Q_{l,k+1}\;-\;Y_{l}(k)]\;.
\end{equation}

The \textit{sending function} is calculated at first by forward recursion, and we substitute $Q_{i,k}=S_{i,k}$ in $(\ref{conserve})$. With this first guess of the amount of traffic in arc i, at time $t_{k+1}$, a first guess of the receiving function can be computed, recursively. 
Finally, we get:

\begin{equation}
C_{i}(k+1)\;=\;F(Y_i(k+1))\;=\;\Psi(Y_i(k),W(k+1))\;,\;\forall\;i\;\in\;\{1,...,n\}.
\end{equation}

\bigskip

\begin{center}
\begin{figure}[h]
\includegraphics[scale=0.7]{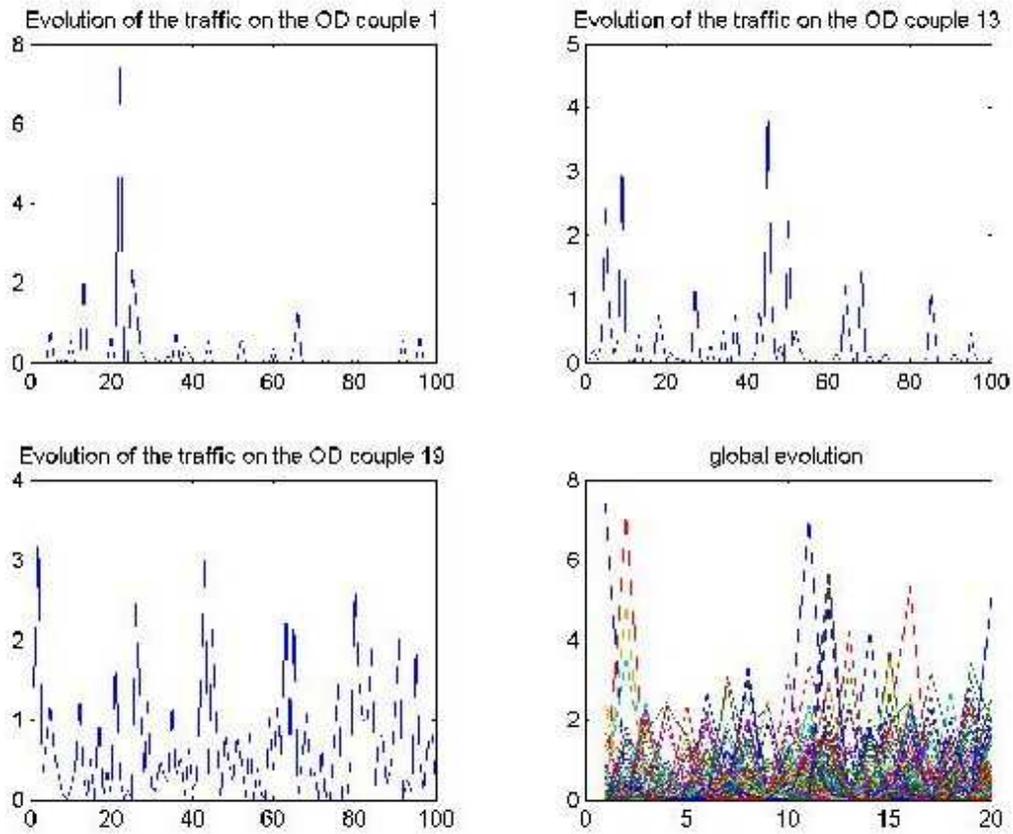}\hfill
\caption{Temporal evolution of the distribution of the traffic on four OD couples.}\label{part1}
\end{figure}

\bigskip

\begin{figure}[h]
\includegraphics[scale=0.7]{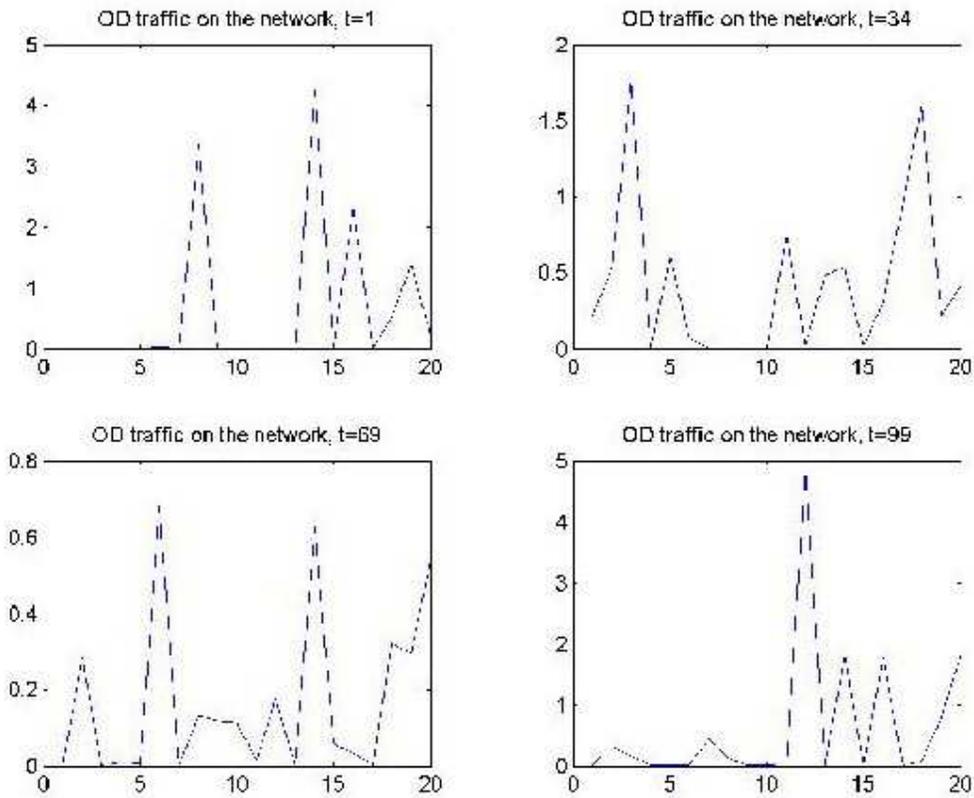}\hfill
\caption{Spatial evolution of the distribution of the traffic at four different instants.}\label{part1}
\end{figure}

\bigskip

\begin{figure}[h]
\includegraphics[scale=0.7]{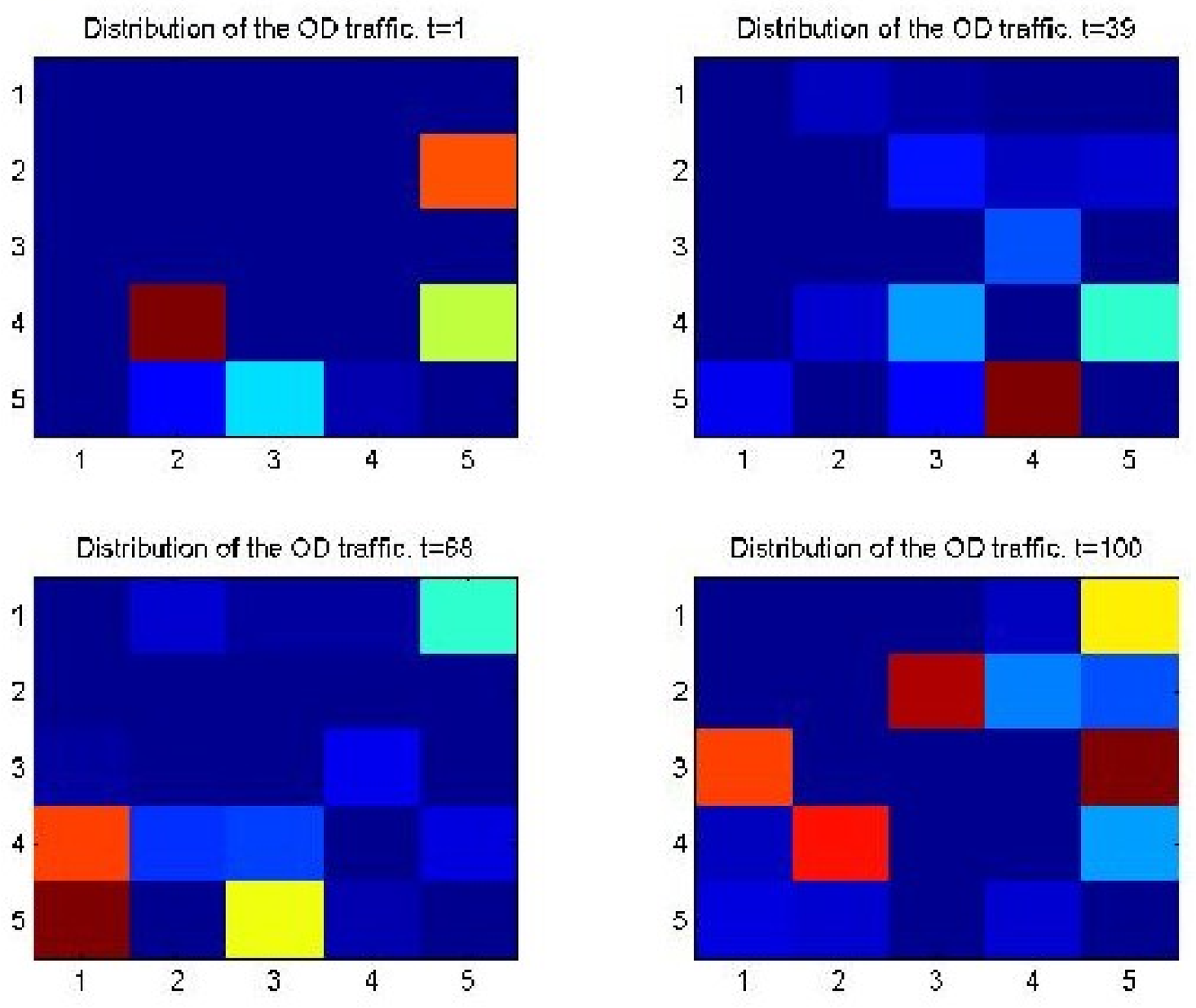}\hfill
\caption{Spatial distribution of the traffic at four instants.}\label{part1}
\end{figure}
\end{center}

\newpage

\section{Conclusion}

We could observe that the performance of the CE method is directly proportionnal to the ratio:
\begin{equation} 
\frac{|\textrm{set of nodes}|}{|\textrm{set of OD couples}|}
\end{equation}

For rather small networks, eg. networks composed of at most $30$ nodes, the CE method works pretty good and suprisingly fastly. 
What's more, it is possible to add some constraints which could guarantee the identifiability of the vector containing the amounts of traffic on the OD couples.
At the end of the estimation part, we get estimators of OD volumes and implicitly, of the routing matrix. In fact, this application is a great proof of the simplicity and versatility of the CE method.\\ 
However, some points remain difficult to tackle.
For example, when the ratio becomes larger than $19$, the CE method performs rather poorly. Furthermore, R. Rubinstein recommand that the sample size of the CE algorithm should be of the form: $$N\;=\;\kappa\;n,\;5\;\leq\;\kappa\;\leq\;10\;.$$
Suppose for example, that that we have to deal with a network of $50$ nodes. Then, at each step of the algorithm we will have to generate a sample of $2450*20000$ vectors. 
Which is completly impossible due to the limited capacities of our computers. But, it is certainly possible to improve the algorithm so as to adapt dynamically the sample size to solve this problem. Nevertheless, the question remains open.
Fortunately, in every network, some specific constraints need to be taken into account. These constraints aim at decreasing the number of unknown parameters. The idea to impose that some OD couples remain equal to zero is an approach, but there are many others. For example, we may want to maximize the global entropy, or other common criteria.\\
Particle Filtering is an efficient and subtle technique to dynamically predict the evolution of the distribution of the traffic on the OD couples for rather small networks.\\
The approaches we use are rather simple to put in application. Nevertheless, they rely on many small parameters which are quite difficult to optimize. Furthermore, the size of the network is still a problem and may be the next challenge of this reflexion.

\newpage

\section{References}
$[1]$ RUBINSTEIN Reuven Y., KROESE Dirk P., \textit{The Cross-Entropy Method}, Springer, $2004$.\\
$[2]$ HU T.C., SHING M.T., \textit{Combinatorial Algorithms}, Dover Publications, second edition, $2002$.\\
$[3]$ SCHRJVER A., \textit{Theory of Linear and Integer Programming}, John Wiley, $2000$.\\
$[4]$ RARDIN R. L., \textit{Optimization in Operations Research}, Prentice Hall, $1998$.\\
$[5]$ DOUCET A., MASKELL S., GORDON N., \textit{Particle Filters for Sequential Bayesian Inference}, Tutorial ISIF, $2002$.\\
$[6]$ MIHAYLOVA L., BOEL R., \textit{A Particle Filter for Freeway Traffic Estimation}.\\
$[7]$ ARULAMPALAM S., MASKELL S., GORDON N., CLAPP T., \textit{A Tutorial on Particle Filters for On-line Non-linear/Non-Gaussian Bayesian Tracking}, IEEE, $2001$.\\
$[8]$ CAMPILLO F., LE GLAND F., \textit{Filtrage Particulaire: quelques exemples "avec les mains" et matlab}, S\'eminaire "filtrage particulaire", CNES, $2003$.\\
\end{document}